\documentclass[twoside]{article}
\usepackage{amsmath}
\usepackage{amssymb}
\usepackage{latexsym}
\newcommand{\mh}{\mathbb}
\newcommand{\mr}{\mathrm}
\newcommand{\mc}{\mathcal}

\newtheorem{thm}{Theorem}

\newtheorem{lem}[thm]{Lemma}
\newtheorem{prop}[thm]{Proposition}

\begin{document}

\begin{center}
\huge Periodic cyclic homology of Hecke algebras and
their Schwartz completions\\[5mm]
\Large Maarten Solleveld\\[1cm]
\end{center}

\begin{minipage}{13cm}
\textbf{Abstract.}
We show that the inclusion of an affine Hecke algebra in its 
Schwartz completion induces an isomorphism on periodic cyclic
homology.
\end{minipage}
\\[1cm]

\noindent
Let $\mc O(V)$ and $C^\infty (X)$ be the algebras of regular 
functions on a nonsingular affine complex variety $V$ and of 
smooth (complex valued) functions on a differentiable manifold 
$X$. The Hochschild-Kostant-Rosenberg theorem \cite{HKR} states 
that there is a natural isomorphism
\begin{equation}\label{eq:1}
HH_* (\mc O(V)) \cong \Omega^* (V)
\end{equation} 
between the Hochschild homology of $\mc O(V)$ and the algebra
of differential forms on $V$, both in the algebraic sense. The 
smooth analogue of this theorem, due to \cite[\S II.6]{Con}, is 
\begin{equation}\label{eq:2}
HH_* (C^\infty(X)) \cong \Omega^* (X; \mh C)
\end{equation}
but now both sides must be interpreted in the topological sense.
\footnote{As concerns the notation, $V$ is a complex algebraic 
variety, so functions and differential forms on $V$ automatically
have complex values. On the other hand, $X$ is a real manifold, and
while it is customary to write $C^\infty (X)$ for $\mh C$-valued
functions, the author believes that it should be mentioned if 
differential forms (and De Rham cohomology) are considered with
complex coefficients.}
Moreover the exterior differential $d$ on $\Omega^*$ corresponds
to the map $B$ on $HH_*$, which implies that 
\begin{align}
\label{eq:3} HP_* (\mc O(V)) &\cong H^*_{DR} (V) \\
\label{eq:4} HP_* (C^\infty (X)) &\cong H^*_{DR} (X; \mh C)
\end{align}
where the right hand sides are $\mh Z / 2 \mh Z$-graded.
However, periodic cyclic homology is much more flexible than 
Hochschild homology and therefore the conditions on $V$ and $X$ 
can be relaxed. In particular \eqref{eq:3} still holds if $V$ is 
singular \cite[Theorem 5]{FT} and \eqref{eq:4} is also valid for 
orbifolds $X$ \cite[\S 4]{Was}.

Suppose now that $X$ is a (smooth) deformation retract of $V$,
endowed with its analytic topology. Because the algebraic and
analytic De Rham cohomologies of $V$ are naturally isomorphic
\cite[Theorem IV.1.1]{Har}, the inclusion $X \to V$ induces
isomorphisms
\begin{align}
\label{eq:5} H^*_{DR}(V) &\to H^*_{DR}(X; \mh C) \\
\label{eq:6} HP_* (\mc O(V)) &\to HP_* (C^\infty (X))
\end{align}
Notice that if $X$ is a compact set of uniqueness for $V$ then 
$C^\infty (X)$ is a completion of $\mc O(V)$. This is remarkable 
since (contrarily to topological $K$-theory) cyclic homology 
theories behave badly with respect to completing algebras. 

For example consider the $C^*$-completion $C(X)$ of $C^\infty(X)$. 
Applying \cite[Section 1]{Joh} to \cite[Corollary 4.9]{Kam} 
we see that
\begin{equation}
HH_n (C(X)) = 0 \qquad \mr{for}\; n > 0
\end{equation}
Hence also $HP_1 (C(X)) = 0$ and
\begin{equation}
HP_0 (C(X)) = HH_0 (C(X)) = C(X)
\end{equation}

In our main theorem we will show that \eqref{eq:6} also holds for
a certain class of noncommutative algebras, namely affine Hecke
algebras and their Schwartz completions. The reader is referred to
the work of Delorme and Opdam \cite{Opd,DO1,DO2} for a precise
definition and a thorough study of the representation theory of
these algebras. One of the first things to notice is that an affine
Hecke algebra is of finite rank over its center, so that we can
use the powerful theory of finite type algebras, which was developed
by Baum, Kazhdan, Nistor and Schneider \cite{BN, KNS}. The author
was particularly inspired by \cite[Theorem 8]{BN}:

\begin{thm}\label{thm:1}
Let $L \to J$ be a spectrum preserving morphism of finite type
algebras. Then the induced map $HP_*(L) \to HP_*(J)$ is an
isomorphism.
\end{thm}

So, just as for commutative finitely generated algebras, the 
periodic cyclic homology of a finite type algebra depends only on 
its spectrum, endowed with Jacobson topology. Unfortunately the 
spectrum $\hat{\mc H}$ of an affine
Hecke algebra $\mc H$ is a rather ugly topological space, it is
a kind of non-separated scheme over $\mh C$. Similarly the spectrum
$\hat{\mc S}$ of the associated Schwartz algebra is a non-Hausdorff
manifold. Notwithstanding these topological inconveniences, it 
follows from \cite{DO2} that we can stratify these spectra so that 
$\hat{\mc S}$ becomes a deformation retract of $\hat{\mc H}$. 
Along these lines we will prove

\begin{thm}\label{thm:2}
Let $\mc H$ be an affine Hecke algebra and $\mc S$ its Schwartz 
completion. Then the inclusion $\mc H \to \mc S$ induces an
isomorphism $HP_* (\mc H) \to HP_* (\mc S)$.
\end{thm}

But first we consider some possible consequences of this theorem.
Let $F$ be a nonarchimedean local field, e.g. a $p$-adic field.
Let $G$ be the group of $F$-rational points of a connected 
reductive algebraic group, and $\mc B(G)$ the set of Bernstein 
components of the smooth dual of $G$ \cite{BD}. The Hecke
algebra $\mc H(G)$ consists of all compactly supported locally 
constant functions on $G$, and it decomposes naturally as
\begin{equation}
\mc H(G) = \bigoplus_{\Omega \in \mc B(G)} \mc H(G)_\Omega
\end{equation}
Similarly $\mc S(G)$ denotes the Schwartz algebra of all rapidly 
decreasing locally constant functions on $G$, which is also an
algebraic direct sum
\begin{equation}
\mc S(G) = \bigoplus_{\Omega \in \mc B(G)} \mc S(G)_\Omega
\end{equation}
It is well known that $\mc H(G)_\Omega$ tends to be Morita equivalent
to the twisted crossed product of a finite group and an affine Hecke
algebra, cf. \cite[Section 5]{ABP}. In particular it has been proved
that, for all Bernstein components of $GL(n,F) ,\: \mc H(G)_\Omega$
is Morita equivalent to a certain affine Hecke algebra $\mc H_\Omega$
\cite{BK}. Moreover in this case $\mc S(G)_\Omega$ is Morita 
equivalent to $\mc S_\Omega$, so Theorem \ref{thm:2} implies 
\cite[Theorem 1]{BHP1}:

\begin{thm}
The inclusion $\mc H(GL(n,F)) \to \mc S(GL(n,F))$ 
induces an isomorphism on periodic cyclic homology.
\end{thm}

More generally, in \cite[Conjecture 8.9]{BHP2} it was conjectured 
that $HP_* (\mc H(G)_\Omega) \to HP_* (\mc S(G)_\Omega)$ is always an
isomorphism. Unfortunately we cannot apply the methods in this paper
to the aforementioned twisted affine Hecke algebras, because not
enough is known about their representation theory. Nevertheless this
conjecture might be proved in another way, in connection with the
Baum-Connes conjecture for $G$, see \cite{Laf}, 
\cite[Proposition 9.4]{BHP2} and \cite[Theorem 12]{Sol}.
\\[2mm]

We recall some of the notations of \cite{Opd}. Let $R_0$ be a
finite, reduced root system with Weyl group $W_0$ and set of 
simple roots $F_0$. Let $\mc R = (X,Y,R_0,R_0^\vee,F_0)$ be a root
datum with affine Weyl group $W = X \rtimes W_0$ and length function 
$l: W \to \mh N$. Pick a label function $q: W \to \mh R^+$, which 
may take different values on nonconjugate simple reflections. 
The affine Hecke algebra $\mc H = \mc H(\mc R,q)$ has a 
$\mh C$-basis $N_w$ in bijection with $W$, and the multiplication 
is defined by
\begin{itemize}
\item $N_v N_w = N_{vw}$ if $l(vw) = l(v) + l(w)$
\item $(N_s + q(s)^{-1/2})(N_s - q(s)^{1/2}) = 0$ for a simple
      reflection $s \in W$
\end{itemize}
The adjoint of $h = \sum_w c_w N_w$ is 
$h^* = \sum_w \overline{c_w} N_{w^{-1}}$.
The Schwartz algebra $\mc S = \mc S(\mc R,q)$ consists of all
(possibly infinite) sums $\sum_{w \in W} c_w N_w$ such that
$w \to |c_w|$ is a rapidly decreasing function, with respect to $l$.
It is a nuclear Fr\'echet *-algebra.

For any $P \subset F_0$ we denote by $\mc H_P = \mc H( \mc R_P,q)$ 
the affine Hecke algebra with root datum 
$\mc R_P = (X_P,Y_P,R_P,R_P^\vee,P)$, where 
\begin{align*}
R_P &= \mh QP \cap R_0 &
R_P^\vee &= \mh Q P^\vee \cap R_0^\vee\\
X_P &= X / X \cap (P^\vee)^\perp & 
Y_P &= Y \cap \mh Q P^\vee
\end{align*}
Furthermore we define 
\[X^P = X / X \cap \mh Q P \quad T^P = \mr{Hom} (X^P, \mh C) \quad
  T_P = \mr{Hom} (X_P, \mh C)
\]
Recall that $T^P$ (and $T_P$ as well of course) decomposes into a
unitary and real split part:
\begin{equation}\label{eq:8}
T^P = \mr{Hom} (X^P, S^1) \times \mr{Hom} (X^P, \mh R^+) =
  T^P_u \times T^P_{rs}
\end{equation}
Let $\Delta_P$ be the set of isomorphism classes of discrete series
of $\mc H_P$, and write $\Delta = \bigcup_{P \subset F_0} \Delta_P$.
Denote by $\Xi$ the analytic variety consisting of all triples
$\xi = (P,t,\delta)$ with $P \subset F_0, \delta \in \Delta_P,
t \in T^P$, and let $\Xi_u$ be the compact submanifold obtained by
restricting to $t \in T^P_u$. For every $\xi \in \Xi$ there exists
a so-called generalized minimal principal series representation 
$\pi(\xi)$ of $\mc H$. Its underlying vector space
$V_\xi = V_{\pi(P,t,\delta)}$ does not depend on $t$, and we let
\[\mc V_\Xi = \bigcup_{(P,\delta) \in \Delta} T^P \times 
  V_{\pi(P,t,\delta)} 
\]
be the corresponding vector bundle over $\Xi$.
Let $\mc W$ be the groupoid, over the power set of $F_0$, with
$\mc W_{PQ} = K_Q \times W(P,Q)$, where $K_Q = T^Q \cap T_Q$ and
\[W(P,Q) = \{w \in W_0 : w(P) = Q\} 
\]
This groupoid acts naturally on $\Xi$ from the left, and for every
$g \in \mc W, \xi = (P,t,\delta) \in \Xi$ there exists an 
intertwiner 
\[\pi(g,\xi) : V_\xi \to V_{g(\xi)}
\]
which is rational in $t$. These intertwiners are unique up to 
scalars and for any choice there exist numbers $c(\delta, g, g')$ 
such that 
\[\pi(g,g',\xi) = c(\delta,g,g') \pi(g, g' \xi) \pi(g, \xi)
\]
In general it is not possible to choose the scalars such that all
the $c(\delta,g,g')$ become 1.

The Langlands classification for $\mc H$ yields 
\cite[Corollary 6.19]{DO2}:

\begin{prop}\label{prop:3}
For every $\pi \in \hat{\mc H}$ there exists a unique association
class $\mc W \xi \in \mc W \setminus \Xi$ such that 
\begin{itemize}
\item $\pi$ is isomorphic to a subquotient of 
      $\pi(\xi) = \pi(P,t,\delta)$
\item $|P|$ is maximal with respect to this property
\end{itemize} 
The resulting map $\hat{\mc H} \to \mc W \setminus \Xi$ is
surjective and finite to one. 
\end{prop}

The orbit $\mc W \xi$ is called the tempered central character of
$\pi$, and $\pi$ extends to $\mc S$ if and only if $\xi \in \Xi_u$.
The intertwiners are unitary on $\Xi_u$, so $\mc W$ also acts on the 
sections of the endomorphism bundle of $\mc V_\Xi$ over $\Xi_u$ by
\[g(f)(\xi) = \pi(g,g^{-1}\xi) f (g^{-1}\xi) \pi(g,g^{-1}\xi)^{-1} 
\]
Now we can formulate \cite[Theorem 4.3]{DO1}:

\begin{thm}\label{thm:4}
The Fourier transform defines an isomorphism of pre-$C^*$-algebras
\[\mc S \to C^\infty \left( \Xi_u;\mr{End}\:\mc V_\Xi \right)^\mc W
\]
\end{thm}

At this point the preparations for the proof of Theorem \ref{thm:2} 
really start. To bring things back to the commutative case we 
construct stratifications of the spectra of $\mc H$ and $\mc S$. 
Choose an increasing chain
\[\emptyset = \Delta_0 \subset \Delta_1 \subset \cdots \subset 
  \Delta_n = \Delta
\]
of $\mc W$-invariant subsets of $\Delta$, with the properties
\begin{itemize}
\item if $(P,\delta) \in \Delta_i$ and $|Q| > |P|$ then 
      $\Delta_Q \subset \Delta_i$
\item the elements of $\Delta_i \setminus \Delta_{i-1}$ form
      exactly one association class for the action of $\mc W$
\end{itemize}
To this correspond two increasing chains of ideals
\begin{align*}
\mc H &= I_0 \supset I_1 \supset \cdots \supset I_n = 0 \\ 
\mc S &= J_0 \supset J_1 \supset \cdots \supset J_n = 0 \\
I_i &= \{ h \in \mc H : \pi(P,t,\delta)(h) = 0 \;\mr{if}\: 
  (P,\delta) \in \Delta_i, t \in T^P \} \\
J_i &= \{ h \in \mc S : \pi(P,t,\delta)(h) = 0 \;\mr{if}\: 
  (P,\delta) \in \Delta_i, t \in T_u^P \}
\end{align*}
For every $i$ pick an element $(P_i,\delta_i) \in \Delta_i 
\setminus \Delta_{i-1}$, let $\mc W_i$ be the stabilizer of
$(P_i,\delta_i)$ in $\mc W$ and write $V_i = 
V_{\pi(P_i,t,\delta_i)}$. Then an immediate consequence of 
Theorem \ref{thm:4} is
\begin{equation}\label{eq:7}
J_{i-1} / J_i \cong C^\infty (T_u^{P_i}; \mr{End}\: V_i)^{\mc W_i}
\end{equation}
while from Proposition \ref{prop:3} we see that the spectrum of 
$I_j/I_i$ corresponds to the inverse image of $\Delta_i \setminus 
\Delta_j$ under the projection $\Xi \to \Delta$. Moreover the 
induced map $\widehat {I_{i-1} / I_i} \to \mc W_i \setminus 
T^{P_i}$ is continuous if we consider $\mc W_i \setminus T^{P_i}$
as an algebraic variety and endow $\widehat {I_{i-1} / I_i}$
with the Jacobson topology. (In fact it is the central character
map for this algebra.)

Recall that the functor $HP_*$ satisfies excision, both in the
algebraic \cite{CQ} and the topological \cite{Cun} setting. 
This means that an extension
\[0 \to I \to A \to A/I \to 0
\]
of algebras gives rise to an exact hexagon
\[\begin{array}{ccccc}
HP_0(I) & \to & HP_0(A) & \to & HP_0(A/I) \\
\uparrow & & & & \downarrow \\
HP_1(A/I) & \leftarrow & HP_1(A) & \leftarrow & HP_1(I)
\end{array}\]
Note however that in the topological case we have to restrict 
ourselves to admissible extensions, i.e. those admitting a 
continuous linear splitting. 

Together with the five lemma this means that in order to prove 
Theorem \ref{thm:2} it suffices to show that each inclusion
\[I_{i-1} / I_i \to J_{i-1} / J_i
\] 
induces an isomorphism on periodic cyclic homology. Therefore we 
zoom in on $J_{i-1} / J_i$. By \cite[Corollary 4.34]{Opd} we can 
extend the action of $\mc W_i$ on $C^\infty (T_u^{P_i}; \mr{End}\; 
V_i)$ to a compact, $\mc W_i$-invariant tubular neighborhood $U$ of 
$T_u^{P_i}$ in $T^{P_i}$. We may assume that $U$ is 
$W_i$-equivariantly diffeomorphic to $T_u^{P_i} \times 
[-1,1]^{\dim T_u^{P_i}}$, and because $[-1,1]$ is compact and 
contractible we can even arrange things so that the extended 
intertwiners $\pi(g,\xi)$ are unitary on all of $U$. It turns out
that we can avoid a lot of technical difficulties by replacing
$J_{i-1} / J_i$ by $C^\infty(U; \mr{End}\: V_i)^{\mc W_i}$.
This is justified by the following result, which is an application
of the techniques developed in \cite{Sol}.

\begin{lem}\label{lem:5}
The inclusion $T_u^{P_i} \to U$ and the Chern character induce 
isomorphisms
\[HP_* (J_{i-1}/J_i) \cong HP_* \left( C^\infty(U; \mr{End}\: 
  V_i)^{\mc W_i} \right) \cong K_* \left(C(U; \mr{End}\: 
  V_i)^{\mc W_i} \right) \otimes_{\mh Z} \mh C
\]
\end{lem}
\emph{Proof.}
The second isomorphism follows directly from the density theorem
for $K$-theory and \cite[Theorem 6]{Sol}.
With the help of \cite{Ill} we pick a $\mc W_i$-equivariant 
triangulation $\Sigma \to T_u^{P_i}$ and we construct a closed cover 
\[\{ V_\sigma : \sigma \;\mr{simplex\:of}\: \Sigma \}
\]
as on \cite[p. 9]{Sol}. Also let $D_\sigma$ be the subset of 
$V_\sigma$ corresponding to the faces of $\sigma$. Using the
projection $p_u : U \to T_u^{P_i}$ we get a closed cover 
\[\{ U_\sigma = \sigma \;\mr{simplex\:of}\: \Sigma \}
\]
of $U$, with 
\[U_\sigma = p_u^{-1} (V_\sigma) \cong 
  V_\sigma \times [-1,1]^{\dim T_u^{P_i}}
\]
According to \cite[p. 10]{Sol} it suffices to show that for any
simplex $\sigma$ we have
\begin{equation}\label{eq:9}
HP_* \left( C_0^\infty(V_\sigma, D_\sigma; \mr{End}\: 
V_i)^{\mc W_\sigma} \right) \cong HP_* \left( C_0^\infty (U_\sigma, 
p_u^{-1}(D_\sigma); \mr{End}\: V_i)^{\mc W_\sigma} \right)
\end{equation}
where $\mc W_\sigma$ is the stabilizer of $\sigma$ in $\mc Wi$.
Well, $V_\sigma \setminus D_\sigma$ is $\mc W_\sigma$-equivariantly
contractible by construction, and it is an equivariant retract of
$U_\sigma \setminus p_u^{-1}(D_\sigma) = p_u^{-1} 
(V_\sigma \setminus D_\sigma)$, so we can apply \cite[Lemma 7]{Sol}.
In this context it says that there exists a finite central extension
$G$ of $\mc W_\sigma$ and a linear representation
\[G \ni g \to u_g \in GL (V_i)
\]
such that the Fr\'echet algebras in \eqref{eq:9} are isomorphic to
\begin{align}
\label{eq:10} C_0^\infty (V_\sigma, D_\sigma; \mr{End}\: V_i)^G \\
\label{eq:11} C_0^\infty (U_\sigma, p_u^{-1}(D_\sigma); 
  \mr{End}\: V_i)^G
\end{align}
The $G$-action is given by
\[g(f)(t) = u_g f(g^{-1} t) u_g^{-1}
\]
where we simply lifted the action of $\mc W_\sigma$ on 
$U_\sigma$ to $G$. 

It is clear that the retraction $U_\sigma \to V_\sigma$ induces a 
diffeotopy equivalence between \eqref{eq:10} and \eqref{eq:11}, so 
it also induces the desired isomorphism \eqref{eq:9}. $\qquad \Box$
\\[2mm]

\emph{Proof of Theorem \ref{thm:2}.}
Consider the finite collection $\mc L$ of all irreducible components
of $(T^{P_i})^g$, as $g$ runs over $\mc W_i$. These are all cosets
of complex subtori of $T^{P_i}$ and they have nonempty 
intersections with $T_u^{P_i}$. Extend this to a collection
$\{ L_j \}_j$ of cosubtori of $T^{P_i}$ by including all irreducible
components of intersections of any number of elements of $\mc L$.
Because the action $\alpha_i$ of $\mc W_i$ on $T^{P_i}$ is algebraic
\[\dim \left( \left( T^{P_i} \right)^g \cap \left( T^{P_i} \right)^w 
  \right) < \max \{ \dim \left( T^{P_i} \right)^g ,\:  
  \dim \left( T^{P_i} \right)^w \}
\]
if $\alpha_i (w) \neq \alpha_i (g)$. Define $\mc W_i$-stable 
submanifolds
\[T_m = \bigcup_{j :\; \dim L_j \leq m} L_j \qquad 
  U_m = T_m \cap U
\]
and construct the following ideals
\begin{align*}
A_m &= \{ h \in I_{i-1} / I_i : \pi(P_i,t,\delta_i)(h) = 0 
  \;\mr{if}\; t \in T_m \} \\
B_m &= C^\infty (U, U_m ; \mr{End}\: V_i)^{\mc W_i} \\
C_m &= C (U, U_m ; \mr{End}\: V_i)^{\mc W_i} 
\end{align*}
Now we have $A_n = B_n = C_n = 0$ for $n \geq \dim T^{P_i}$ and
\[A_n = I_{i-1} / I_i \quad 
  B_n = C^\infty (U; \mr{End}\: V_i)^{\mc W_i} \quad 
  C_n = C (U; \mr{End}\: V_i)^{\mc W_i} \quad \mr{for}\; n<0
\]
Using excision and Lemma \ref{lem:5}, it will be sufficient to 
show that the inclusions
\[A_{m-1} / A_m \to B_{m-1} / B_m
\] 
induce isomorphisms on $HP_*$, so let us compute the periodic 
cyclic homologies of these quotient algebras.

Because $T_m$ is an algebraic subvariety of $T^{P_i}$ the spectrum 
of $A_{m-1} / A_m$ consists precisely of the 
irreducible representations of $I_{i-1} / I_i$ with tempered 
central character in $(P_i, T_m \setminus T_{m-1}, \delta_i)$.
We let $r_i(t)$ be the number of $\pi \in \widehat{I_{i-1} / I_i}$
corresponding to $(P_i,t,\delta_i)$. 
From the proof of \cite[Proposition 6.17]{DO2} and we see that 
$r_i(t |t|^s) = r_i(t) \; \forall s > -1$, and 
$r_i(t |t|^{-1}) = r_i(t)$ if the stabilizers in $\mc W_i$ of $t$
and $t |t|^{-1}$ are equal. Choose a minimal subset  
$\{ L_{m,k} \}_k$ of $\mc L$ such that every $m$-dimensional element
of $\mc L$ is conjugate under $\mc W_i$ to a $L_{m,k}$. Let 
$\mc W_{m,k}$ be the stabilizer of $L_{m,k}$ in $\mc W_i$ and write
$r_k = r_i(t)$ for some $t \in L_{m,k} \setminus T_u^{P_i}$. Then
the spectrum of $A_{m-1} / A_m$ is homeomorphic to
\begin{align*}
&\bigsqcup_k \bigsqcup_{l=1}^{r_k} \left( L_{m,k} \setminus T_{m-1} 
  \right) / \mc W_{m,k} \qquad = \\
&\bigsqcup_k \bigsqcup_{l=1}^{r_k} \left( L_{m,k} / \mc W_{m,k} 
 \right) \setminus \left( (L_{m,k} \cap T_{m-1}) / \mc W_{m,k} \right)
\end{align*}
Let us call the left hand side of this expression $X_m$, and its
subset on the right hand side $Y_m$; these are complex 
algebraic varieties. By \cite[Proposition 6.17.iii]{DO2} we have
\begin{equation}\label{eq:12}
Z(A_{m-1} / A_m) \cong \{ f \in \mc O(X_m) : f(Y_m) = 0 \}
\end{equation}
so with \cite[Theorem 9]{KNS} we get
\begin{equation}\label{eq:13}
HP_* (Z(A_{m-1} / A_m)) \cong \check H^* (X_m,Y_m;\mh C)
\end{equation}
Now, as $Z(A_{m-1} / A_m) \to A_{m-1} / A_m$ is a spectrum
preserving morphism of finite type algebras, \cite[Theorem 8]{BN}
tells us that
\begin{equation}\label{eq:14}
HP_* (A_{m-1} / A_m) \cong \check H^* (X_m,Y_m;\mh C)
\end{equation}

On the other hand, by \cite[Th\'eor\`eme IX.4.3]{Tou} the extension
\[0 \to C^\infty (U, U_m ; \mr{End}\: V_i) \to C^\infty 
  (U; \mr{End}\: V_i) \to C^\infty (U_m ; \mr{End}\: V_i) \to 0
\]
is admissible, and since $\mc W_i$ is finite the same holds for
\[0 \to B_m \to C^\infty (U; \mr{End}\: V_i)^{\mc W_i} \to 
  C^\infty (U_m ; \mr{End}\: V_i)^{\mc W_i} \to 0
\]
So by \cite[Theorem 6]{Sol} we have isomorphisms
\begin{equation}
HP_* (B_m) \leftarrow K_*(B_m) \otimes_{\mh Z} \mh C \to 
K_*(C_m) \otimes_{\mh Z} \mh C
\end{equation}
By construction the stabilizer in $\mc W_i$ of $t \in U$ is constant
on the connected components of $U_m \setminus U_{m-1}$, and, by
the continuity of the intertwiners $\pi(g,P_i,t,\delta_i)$, the same
can be  said of the type of $V_i$ as a 
$\left< g \right>$-representation (on $U_m^g$). 
Thus, with $L_{m,k}$ as above, we get
\begin{align*}
C_{m-1} / C_m &\cong \bigoplus_k C_0 \left( (L_{m,k} \cap U_m) 
\setminus (L_{m,k} \cap U_{m-1}) / \mc W_{m,k} \right) \otimes S_k \\
&= \bigoplus_k C_0 \left( (L_{m,k} \cap U_m) / \mc W_{m,k}, 
 (L_{m,k} \cap U_{m-1}) / \mc W_{m,k} \right) \otimes S_k 
\end{align*}
where the $S_k$ are certain finite dimensional semisimple 
$\mh C$-algebras. Because $I_{i-1} / I_i$ is dense in 
$C^\infty (U; \mr{End}\: V_i)^{\mc W_i}$ we must have 
$\dim Z(S_k) = r_k$. Consequently
\begin{equation}\label{eq:15}\begin{split}
HP_* (B_{m-1}/B_m) &\cong K_* (C_{m-1}/C_m) \otimes_{\mh Z} \mh C \\
&\cong \bigoplus_k \check H^* \left((L_{m,k} \cap U_m) / \mc W_{m,k}, 
  (L_{m,k} \cap U_{m-1}) / \mc W_{m,k} ; \mh C \right)^{r_k} \\
&= \check H^* \left( \bigsqcup_k \bigsqcup_{l=1}^{r_k} (L_{m,k} 
  \cap U_m) / \mc W_{m,k}, \bigsqcup_k \bigsqcup_{l=1}^{r_k}(L_{m,k} 
  \cap U_{m-1}) / \mc W_{m,k} ; \mh C \right) \\
&= \check H^* (X'_m ,Y'_m ; \mh C)
\end{split}\end{equation}
where $X'_m := X_m \cap U / \mc W_i$ and 
$Y'_m := Y_m \cap U / \mc W_i$.

It follows from this and the density theorem that the inclusions
\begin{equation}\label{eq:16}
C_0^\infty (X'_m ,Y'_m ) \to C_0 (X'_m ,Y'_m ) \cong
Z (C_{m-1} / C_m ) \to C_{m-1} / C_m
\end{equation}
induce isomorphisms on $K$-theory with complex coefficients. 
From \eqref{eq:13} - \eqref{eq:16} we construct the commutative 
diagram
\begin{equation}
\begin{array}{ccccc}
HP_* (A_{m-1} / A_m) & \leftarrow & HP_* \big( Z(A_{m-1} / 
 A_m) \big) & \to & \check H^* (X_m,Y_m;\mh C) \\
\downarrow & & \downarrow & & \downarrow \\
HP_* (B_{m-1} / B_m) & \cong & HP_* \big( C_0^\infty 
 (X'_m ,Y'_m ) \big) & \to & \check H^* (X'_m , Y'_m; \mh C)
\end{array}
\end{equation}
The pair $(X'_m, Y'_m)$ is a deformation retract of $(X_m,Y_m)$, 
so all maps in this diagram are isomorphisms. 
Working our way back up, using excision, we find that also
\[HP_* (I_{i-1}/I_i) \to HP_* \left( C^\infty (U; \mr{End}\: 
  V_i)^{\mc W_i} \right) \to HP_* (J_{i-1}/J_i)
\]
and finally 
\[HP_* (\mc H) \to HP_*(\mc S)
\] 
are isomorphisms. $\qquad \Box$
\\[1cm]

\textbf{Acknowledgements.}
The author would like to thank Eric Opdam for many illuminating
conversations on this subject.
\\[8mm]

\begin{minipage}{13cm}
Korteweg-de Vries Institute for Mathematics\\
Universiteit van Amsterdam\\
Plantage Muidergracht 24\\
1018TV  Amsterdam\\
The Netherlands\\
Email: mslveld@science.uva.nl
\end{minipage}

\end{document}